\documentclass[%
 aip,
cp,  % Conference Proceedings
 amsmath,amssymb,%nobibnotes,
 reprint,%
]{revtex4-2}

\usepackage{graphicx}% Include figure files
\usepackage{dcolumn}% Align table columns on decimal point
\usepackage{bm}% bold math
\usepackage[utf8]{inputenc}
\usepackage[T1]{fontenc}
    \usepackage{mathptmx}
    \usepackage{amsmath}
    \usepackage{mathtools}  %for double arrow

    \graphicspath{ {./Figures/} } 
    \usepackage{lineno}
\begin{document}

\title{
    Solution of Linear Systems of Equations \textbf{Ax=b} and \textbf{Ax=0} using Unifying Approach with Geometric Algebra: \\ 
    Outer Product Application and Angular Conditionality
    }% Force line breaks with \\

\author{Vaclav Skala} % Write as First name Surname
 \email[Corresponding author: ]{skala@kiv.zcu.cz \quad www.VaclavSkala.eu }
%\author{Author's Name}%
% \email{second.author@institution.edu.}
\affiliation
{Department of Computer Science and Engineering \\
    Faculty of Applied Sciences, 
    University of West Bohemia \\
    Univerzitni 8
    CZ 306 14 Pilsen Czech Republic
  % Replace this text with an author's affiliation
  % (use complete addresses, including country name or code).% Force line breaks with \\ if necessary
}

%\author{Another's Name}
% \email{third.author@anotherinstitution.edu}
%\affiliation{Second institution and/or address% Force line breaks with \\ if necessary }%
%\affiliation{You would list an author's second affiliation (if applicable) here.}

% \date{\today} % It is always \today, today, but any date may be explicitly specified
              % Not printed for conference proceedings

\begin{abstract}
    A solution of linear systems of equations \textbf{Ax=b} and \textbf{Ax=0} is a vital part of many computational packages.
    This paper presents a novel formulation based on the projective extension of the Euclidean space using the outer product (extended cross-product). 
    This approach enables to solve the both cases, i.e. 
    \textbf{Ax=b} and \textbf{Ax=0}
    The proposed approach leads actually to an "analytical" solution of linear systems in the form 
    $\bm{\xi}=\bm{\alpha}_1 \wedge \bm{\alpha}_2 \wedge \ldots \wedge \bm{\alpha}_n,$ 
    on which the other vector operation can be applied before using the numerical evaluation.
    
    This contribution also proposes a new approach to the conditionality estimation of matrices applicable to non-squared matrices.
    It splits the conditionality to "structural" conditionality showing matrix property if nearly unlimited precision is used, "numerical" issue which depends on numerical representation with respect to the right-hand side influence, if given.
\end{abstract}


\maketitle

\section{Introduction \protect}
    Solutions of a linear system of equations is a vital part of a solution of many computational packages.
    There are two types of linear systems 
    \textbf{Ax=b} and \textbf{Ax=0}, i.e. with the right-hand side and without it. 
    In the first case, the matrix $\mathbf{A}$ of the size $n \times n$ is expected to be non-singular, i.e. $\det(\mathbf{A}) \ne 0$, and the matrix is to be positive definite, if an iterative solver is to be used.
    This is related to the explicit formulations.
    In the second case, which is related to the implicit formulations, the matrix  $\mathbf{A}$ of the size $n \times (n+1)$.
    It leads to one-dimensional parametric solution.
    
    Methods of the linear system of equations $\mathbf{Ax=b}$ solutions have been deeply studied and very sophisticated methods have been developed. 
    However, the limited precision of real numbers representation using the IEEE-754 standard \cite{wiki:IEEE-754}\cite{IEEE754-219} leads to severe problems with numerical stability, robustness, speed of computation and even with the correctness of the solution, especially with the growing size $n$ of the matrix $\mathbf{A}$. 
    It should be noted that the size $n$ can be quite high, e.g. $10^6$ and higher, see Majdisova\cite{Majdisova201751}.
    
    However, many engineering problems lead to the "ill-conditioned" matrices. The conditionality $\kappa(\mathbf{A})$ of a matrix $\mathbf{A}$ can be  estimated as $\kappa(A) = | \lambda_{max}|/|\lambda_{max}|$, where $\lambda_i \in C^1$ are eigenvalues of the matrix $\mathbf{A}$, which might be evaluated, e.g. using the Gershgorin circle theorem\cite{enwiki:Gershgorin-2021}.
    A typical example of the very ill-conditioned matrix is the Hilbert matrix\cite{enwiki:Hilbert-2021}.
    
    Methods of the linear system of equations $\mathbf{Ax=0}$ solutions are partially out of the main research interest. 
    In this case, when the matrix $\mathbf{A}$ is of the size $n \times (n+1)$, the linear system $\mathbf{Ax=0}$ represents a solution of many physical problems, seemingly in a one parametric form, which is difficult to formulate analytically.
    The conditional issues are not quite well defined and analyzed in this case 
    (mostly only linear independence of rows is evaluated regardless of the numerical precision available).

% ======================================      
    
\section{Projective extension of the Euclidean space \protect}  
    The concept of the projective extension of the Euclidean space, i.e. the projective space, was originated from the visual perception of parallel lines which seem to meet in infinity.  
    It uses homogeneous coordinates and two equivalent forms can be found:
    \begin{itemize}
        \item the form $[x_1,\ldots,x_n:x_w]$ is mostly used in the computer graphics related fields, namely 
        $[x,y:w]$ in the case of $P^2$, resp. $[x,y,z:w]$ in the case of $P^3$, where $w$ is the homogeneous coordinate.
        \item the form $[x_0:x_1,\ldots,x_n]$ is used in the mathematical fields and the $x_0$ is the homogeneous coordinate. This form has the advantage that the homogeneous coordinate is on the first position.
    \end{itemize}
    It should be noted that "$:$" is used to emphasize that the $x_w$, resp $x_0$ has a different meaning as it is actually the "scaling factor", i.e. without a physical unit, while $x_1,\ldots,x_n$ has different physical units, e.g. meters[m] etc.
    
\begin{figure}[!ht]
        \centering
        \includegraphics[scale=0.4]{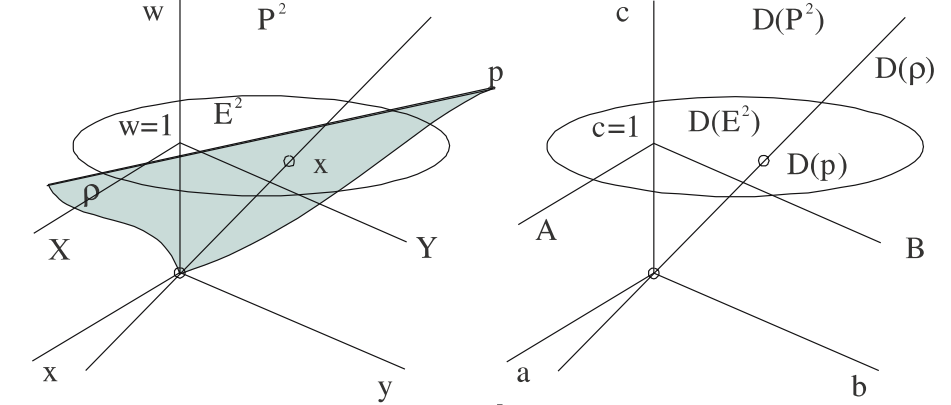}
        \caption{The Euclidean space, the projective space and the dual representation}
        \label{fig:Projective-space}
    \end{figure}

    The mutual conversion between the Euclidean space and projective space is given as:
    \begin{equation}
    \begin{split}
        X_i = \frac{x_i}{x_0} \quad x_0 \ne 0 \quad ,\ resp. \quad X_i = \frac{x_i}{x_w} \quad 
        x_w \ne 0 \quad , \quad i=1,\ldots,n
    \end{split}
    \end{equation}
    where $X_i$ are coordinates in the Euclidean space. \\
    In the case of the $E^2$ space
    \begin{equation}
    \begin{split}
        X = \frac{x}{x_0} \quad Y = \frac{Y}{x_0} \quad x_0 \ne 0  \quad  ,\  resp. \\
        \quad X = \frac{x}{w} \quad \quad Y = \frac{y}{w} \quad w \ne 0 
    \end{split}
    \end{equation}    
    where $(X,Y)$, resp.$[x,y:w]$ are coordinates in the Euclidean space $E^2$, resp.in the projective space $P^2$. The extension to the $E^3$, resp. $E^n$ space is straightforward, see Vince\cite{Vince-2010}\cite{Yamaguchi-2002}. 
    \\
    The geometrical interpretation of the Euclidean  ($x_w=1$, resp. $x_0=1$) and the projective spaces is presented at Fig.\ref{fig:Projective-space}.

\section{Geometric, Inner and Outer Products \protect}
    
    The \textit{Geometric Algebra} (GA)  introduces a "new" product called \textit{geometric product}, 
    which composes the \textit{dot product} and \textit{outer product} as follows:
    \begin{equation}\label{geometric-product}
    \begin{split}
        \mathbf{ab}= \mathbf{a} \cdot \mathbf{b} + \mathbf{a} \wedge \mathbf{b}
    \end{split}
    \end{equation}
    where $\mathbf{ab}$ is the new entity, $\mathbf{a}$, $\mathbf{b}$ are "movable" vectors (in the mathematical sense) in the $E^n$ space, "$\cdot$" means the \textit{dot product} and "$\wedge$" means the \textit{outer product}, see Vince\cite{Vince-2008}\cite{Vince-2009}.

    It is a "set of objects" with different dimensionalities and properties, in general. 
    In the case of the $n$-dimensional space, the vectors are defined as $\mathbf{a}=(a_1 \mathbf{e}_1+...+a_n \mathbf{e}_n)$, $\mathbf{b}=(b_1 \mathbf{e}_1+...+b_n \mathbf{e}_n)$, 
    where the $\mathbf{e}_i$ vectors form orthonormal vector basis in $E^n$. 
    In the $E^3$ case, the following objects can be used in geometric algebra Vince\cite{Vince-GA}, Macdonald\cite{Macdonald2017853}, Doran\cite{Doran-2002}, Dorst\cite{Dorst-2011}, Katani\cite{Kanatani20151}, Hildebrand\cite{Hildebrand-2013}
    :
    \begin{center}
    \begin{tabular}{ c c c c c }
     1\quad & 0-vector (scalar) & \hspace{1cm} & $\mathbf{e}_{12}$, $\mathbf{e}_{23}$, $\mathbf{e}_{31}$  & 2-vectors (bivectors)\\ 
     $\mathbf{e}_1, \mathbf{e}_2, \mathbf{e}_3$ \quad & 1-vector (vectors) \quad & \hspace{1cm} & $\mathbf{e}_{123}$ & 3-vector (pseudoscalar)
    \end{tabular}
    \end{center}
    The significant advantage of the geometric algebra is, that it is more general than the Gibbs algebra as it can handle all objects with dimensionality up to $n$. 
    The geometry algebra uses the following operations, including the inverse of a vector.
    \begin{equation}
         \mathbf{a} \cdot \mathbf{b} = \frac{1}{2}\mathbf{(ab+ba)} \hspace{1cm}
         \mathbf{a}\wedge \mathbf{b} = -\mathbf{b}\wedge \mathbf{a} \\ \hspace{1cm} \mathbf{a}^{-1}=\mathbf{a}/||\mathbf{a}||^2
    \end{equation}
     % ----
    It should be noted, that the geometric algebra is \textit{anti-commutative} and the “pseudoscalar”  $I$  has the basis $\mathbf{e}_1 \mathbf{e}_2 \mathbf{e}_3$ (briefly as $\mathbf{e}_{123}$) in the $E^3$case, i.e.
    \begin{equation}
        \mathbf{e}_i \mathbf{e}_j = -\mathbf{e}_j \mathbf{e}_i \hspace{1cm} \mathbf{e}_i \mathbf{e}_i = 1 \hspace{1cm} \mathbf{e}_1 \mathbf{e}_2 \mathbf{e}_3 = I \hspace{1cm} 
        \mathbf{a} \wedge \mathbf{b} \wedge \mathbf{c} = q
    \end{equation}
    where $q$ is a scalar value (actually a pseudoscalar).
    In the case of the $E^3$ case, the equation Eq.\ref{geometric-product} is equivalent to:
    \begin{equation}
    \begin{split}
        \mathbf{ab}= \mathbf{a} \cdot \mathbf{b} + \mathbf{a} \times \mathbf{b}
    \end{split}
    \end{equation}
    where "$\times$" mean the cross-product.
    \\   
    In general, the geometric product is represented as:
    \begin{equation}
        \mathbf{ab}=\sum_{i,j=1}^{n,n}a_{i}\mathbf{e}_{i}b_{j} \mathbf{e}_{j} \hspace{1cm} 
        \mathbf{a} \cdot \mathbf{b}=\sum_{i=1}^{n,n}a_{i}\mathbf{e}_{i}b_{i}\mathbf{e}_{i} 
    \end{equation}
    \begin{equation}
        \mathbf{a} \wedge \mathbf{b} = \sum_{i,j=1 \& i\neq j}^{n,n}a_{i}\mathbf{e}_{i}b_{j}\mathbf{e}_{j} = \sum_{i,j=1,\& i>j}^{n} (a_{i}b_{j}-a_{j}b_{i})\mathbf{e}_{i}\mathbf{e}_{j}
    \end{equation}

    It is not a “user-friendly” notation for practical applications and causes problems in practical implementations, as the geometric product is anti-commutative.
    \\
    The efficient computation of the geometric product $\mathbf{ab}$ of two vectors $\mathbf{a}$ and $\mathbf{b}$ using the tensor product WiKi\cite{enwiki:Tensor-2021} defined by Eq.\ref{Eq:Tensor-product} was described in Skala\cite{Skala2022437}
    \begin{equation}\label{Eq:Tensor-product}
    \begin{split}
        \mathbf{a} \otimes \mathbf{b} = 
        \begin{bmatrix}
            a_{1}b_{1} & a_{1}b_{2} & \cdots  & a_{1}b_{m}  \\
            a_{2}b_{1} & a_{2}b_{2} & \cdots  & a_{2}b_{m}  \\
            \vdots     &    \vdots  & \ddots  & \vdots   \\
            a_{n}b_{1} & a_{n}b_{2} & \ldots  & a_{n}b_{m}
        \end{bmatrix}
    \end{split}    
    \end{equation}
In the case of the $E^3$ space, it should be noted that the matrix $\mathbf{Q}$ has the following combinations of the basis vectors:
    \begin{equation}
    \begin{split}
            \mathbf{Q} = 
            \begin{bmatrix}
                \mathbf{e}_1 \mathbf{e}_1 & \mathbf{e}_1 \mathbf{e}_2 & \mathbf{e}_1 \mathbf{e}_3 \\
                \mathbf{e}_2 \mathbf{e}_1 & \mathbf{e}_2 \mathbf{e}_2 & \mathbf{e}_2 \mathbf{e}_3 \\
                \mathbf{e}_3 \mathbf{e}_1 & \mathbf{e}_3 \mathbf{e}_2 & \mathbf{e}_3 \mathbf{e}_3 \\
            \end{bmatrix} =
            \begin{bmatrix}
                1 & \mathbf{e}_1 \mathbf{e}_2 & -\mathbf{e}_3 \mathbf{e}_1 \\
                -\mathbf{e}_1 \mathbf{e}_2 & 1 & \mathbf{e}_2 \mathbf{e}_3 \\
                \mathbf{e}_3 \mathbf{e}_1 & -\mathbf{e}_2 \mathbf{e}_3 & 1 \\
            \end{bmatrix}      
    \end{split}
    \end{equation}
    In the $E^3$ case, the right-handed coordinate system has the orthonormal basis
    $\mathbf{e}_1 \mathbf{e}_2$, $\mathbf{e}_2 \mathbf{e}_3$,  $\mathbf{e}_3 \mathbf{e}_1 $  and therefore the value of $q_{13}$ results into the $-\mathbf{e}_3 \mathbf{e}_1$ value.
    
    It means, that the results of the $\mathbf{a} \otimes \mathbf{b}$ operations is: 
    \begin{equation}
    \begin{split}
        \mathbf{a} \otimes \mathbf{b} =
        \begin{bmatrix}
            a_1 b_1 \mathbf{e}_1 \mathbf{e}_1 & a_1 b_2 \mathbf{e}_1 \mathbf{e}_2 & - a_1 b_3\mathbf{e}_3 \mathbf{e}_1 \\
            - a_2 b_1\mathbf{e}_1 \mathbf{e}_2 & a_2 b_2 \mathbf{e}_2 \mathbf{e}_2  & a_2 b_3 \mathbf{e}_2 \mathbf{e}_3 \\
            a_3 b_1 \mathbf{e}_3 \mathbf{e}_1 & - a_3 b_2 \mathbf{e}_2 \mathbf{e}_3 & a_3 b_3 \mathbf{e}_3 \mathbf{e}_3 \\
        \end{bmatrix}      
    \end{split}
    \end{equation}
    including the right-hand orientation of the coordinate system, resulting into the "-" sign in the matrix. 
    \\
    Note, that $\mathbf{e}_i \mathbf{e}_i=1$ by definition and therefore: 
    \begin{equation}
    \begin{split}
        \mathbf{a} \otimes \mathbf{b} =
        \begin{bmatrix}
            a_1 b_1  & a_1 b_2 \mathbf{e}_1 \mathbf{e}_2 & - a_1 b_3\mathbf{e}_3 \mathbf{e}_1 \\
            - a_2 b_1\mathbf{e}_1 \mathbf{e}_2 & a_2 b_2   & a_2 b_3 \mathbf{e}_2 \mathbf{e}_3 \\
            a_3 b_1 \mathbf{e}_3 \mathbf{e}_1 & - a_3 b_2 \mathbf{e}_2 \mathbf{e}_3 & a_3 b_3  \\
        \end{bmatrix}      
    \end{split}
    \end{equation}
    It can be seen, that the diagonal represents the \textit{inner product}, while non-diagonal elements are related to the \textit{outer product}, see the Appendix.    
    \\
    Let us consider the projective extension of the Euclidean space and the use of the homogeneous coordinates. 
    The geometric algebra concept can be extended for the $P^n$ projective space as:
    \begin{equation}
    \begin{split}
        [\mathbf{a}:w_a] [\mathbf{b}:w_b]= [\mathbf{a} \cdot \mathbf{b}:w_a w_b] + [\mathbf{a} \wedge \mathbf{b}:w_a w_b]
    \end{split}
    \end{equation}
    In this case, the values of $\mathbf{a}$ and $\mathbf{b}$ represent some physical entity, e.g. a position in the n-dimensional Cartesian space. 
    It means, that $[\mathbf{a}:w_a]$ and $[\mathbf{b}:w_b]$ are not movable vectors, but they are fixed to the origin of the Cartesian coordinate system.
    \\
    Then the vectors $\mathbf{a}=[a_1,a_2,a_3:a_4 ]^T$ and $\mathbf{b}=[b_1,b_2,b_3:b_4 ]^T$ in the projective space actually represent vectors $(a_1/a_4,a_2/a_4,a_3/a_4)$ and $(b_1/b_4,b_2/b_4,b_3/b_4)$ in the Euclidean space $E^3$.
    \\
    In the following, the homogeneous coordinates will be used.
    Then the geometric product is represented by the tensor product for the 4-dimensional case, see Eq.\ref{Eq:Tensor-product-proj}, as:
    \begin{equation}\label{Eq:Tensor-product-proj}
    \begin{split}
        \mathbf{ab} \xLeftrightarrow[\text{repr}]{\text{}} %arrow symbol
        \mathbf{ab}^T =
        \mathbf{a} \otimes \mathbf{b} = 
        \mathbf{Q} 
        \\
        =
        \begin{bmatrix}
            a_{1}b_{1} & a_{1}b_{2} & a_{1}b_{3} & a_{1}b_{4}\\
            a_{2}b_{1} & a_{2}b_{2} & a_{2}b_{3} & a_{2}b_{4}\\
            a_{3}b_{1} & a_{3}b_{2} & a_{3}b_{3} & a_{3}b_{4} \\
            a_{4}b_{1} & a_{4}b_{2} & a_{4}b_{3} & a_{4}b_{4}   
        \end{bmatrix} 
        = 
        \begin{bmatrix}
            a_{1}b_{1} & a_{1}b_{2} & a_{1}b_{3} & \vline & a_{1}w_b\\
            a_{2}b_{1} & a_{2}b_{2} & a_{2}b_{3} & \vline & a_{2}w_b\\
            a_{3}b_{1} & a_{3}b_{2} & a_{3}b_{3} & \vline & a_{3}w_b \\
            \hline 
             w_a b_{1} &  w_a b_{2} &  w_a b_{3} & \vline & w_a w_b 
        \end{bmatrix} 
%        \\
        =
        \begin{bmatrix} 
            &            & \vline &    \\
            & \mathbf{P} & \vline & \mathbf{r} \quad &\\
            \hline 
            &           & \vline \\
            & \mathbf{s}^T  & \vline &  w \\
        \end{bmatrix} 
    \end{split}
    \end{equation}
    where  $\mathbf{P} = \mathbf{B+U+D} $ are $B$ottom triangular, $U$pper triangular, $D$iagonal matrices, $a_4,b_4$ are the homogeneous coordinates, i.e. actually $w_a, w_b$ (will be explained later), and the operator $\otimes$ means the anti-commutative tensor product.

    Let us consider the projective extension of the Euclidean space and use of the homogeneous coordinates. Let us consider vectors $\mathbf{a}=[a_1,a_2,a_3:a_4 ]^T$ and $\mathbf{b}=[b_1,b_2,b_3:b_4 ]^T$, which represents actually vectors $(a_1/a_4,a_2/a_4,a_3/a_4)$ and $(b_1/b_4,b_2/b_4,b_3/b_4)$ in the $E^3$ space. It can be seen, that the diagonal $\mathbf{D}$ of the matrix $\textbf{P}$ actually represents the inner product in the projective representation:
    \begin{equation}
    \begin{split}
        \mathbf{a} \cdot \mathbf{b} = [(a_1 b_1 + a_2 b_2 + a_3 b_3):a_4 b_4]^T \\
        \triangleq \frac{a_1 b_1 + a_2 b_2 + a_3 b_3}{a_4 b_4} 
        = \frac{a_1 b_1 + a_2 b_2 + a_3 b_3}{w_a w_b}
    \end{split}
    \end{equation}
    where $\triangleq$ means projectively equivalent. The inner product actually represents trace $tr(\mathbf{P})$ of the matrix $\mathbf{Q}$.
    
    The outer product (the cross-product in the $E^3$ case) is then represented 
    by matrices $\mathbf{B} + \mathbf{D}$ 
    respecting anti-commutativity as:
    \begin{equation}
    \begin{split}
        \mathbf{a} \wedge \mathbf{b} \xLeftrightarrow[\text{repr}]{\text{}} [\sum_{i,j=1\&i>j}^{3,3} (a_i b_j \mathbf{e}_i \mathbf{e}_j - b_i a_j \mathbf{e}_i \mathbf{e}_j ) :a_4 b_4]^T
        \\
        \triangleq \frac{\sum_{i,j \& i>j}^{3,3} (a_{i}b_{j} - b_{i}a_{j}) \mathbf{e}_{i} \mathbf{e}_{j}}{a_4 b_4} 
        = \frac{\sum_{i,j \& i>j}^{3,3} (a_{i}b_{j} - b_{i}a_{j}) \mathbf{e}_{i} \mathbf{e}_{j}}{a_w b_w}
    \end{split}
    \end{equation}
    ??????\\
    It should be noted, that the outer product can  be used for a solution of a linear system of equations  $\mathbf{Ax=b}$ or $\mathbf{Ax=0}$, too.
    \\ ???\\\\ 
    
    \begin{equation}
    \begin{split}
        \mathbf{Q} - \mathbf{Q}^T = 
        \begin{bmatrix}
                      0        & a_1 b_2 - a_2 b_1 & a_1 b_3 - a_3 b_1 & \vline & a_1 w_b - w_a b_1\\
            a_2 b_1 - a_1 b_2  &         0         & a_2 b_3 - a_3 b_2 & \vline & a_2 w_b - w_a b_2\\
            a_3 b_1 - a_1 b_3  & a_3 b_2 - a_2 b_3 &        0          & \vline & a_3 w_b - w_a b_3\\
            \hline 
             w_a b_1 - a_1 w_b & w_a b_2 - a_2 w_b &  w_a b_3 - a_3 w_b & \vline & 0
        \end{bmatrix} 
        = 
        \begin{bmatrix} 
            &            & \vline &    \\
            & \mathbf{P} & \vline & \mathbf{r} \quad &\\
            \hline 
            &           & \vline \\
            & \mathbf{s}^T  & \vline &  w \\
        \end{bmatrix}
    \end{split}
    \end{equation}
    where $\mathbf{P} = \mathbf{B} + \mathbf{U} + \mathbf{D} $
    It can be seen, that the matrix $\mathbf{U}$ represents the outer product, multiplied by $w = w_a w_b$.
    The vector $\mathbf{r}$ represents difference 
    $[\mathbf{a}:w_a] - [\mathbf{b}:w_b]$ and the vector
    $\mathbf{s}$ represents the difference 
    $ [\mathbf{b}:w_b]- [\mathbf{a}:w_a]$, actually multiplied by $w = w_a w_b$.
    It means, that the vector $\mathbf{r}$ resp. $\mathbf{s}$ represents a directional vector of a line passing two points in $E^2$ space. 
    Also, a relationship with the P\:ucker coordinates can be seen.

\section{Principle of Duality \protect}
    The  duality principle is very important principle, but unfortunately usually not covered in introductory mathematical courses.
    The principle of duality is important principle, in general. 
    The projective principle of duality states that any theorem remains true when we interchange the words “point” and “line” in the $P^2$ , resp. “point” and “plane” in $P^3$, “lie on” and “pass-through”, “join” and “intersection” and so on. 
    Once the theorem has been established, the dual theorem is obtained as described above Johnson\cite{Johnson-1996}\cite{enwiki:Duality}. 
    In other words, the principle of duality in $P^2$ says that in all theorems, it is possible to substitute a term “point” by a term “line” and term “line” by the term “point” and the given theorem stay valid; similarly in $P^3$ space with term "point" and "plane".
    \\ 
    \\
    Applying the principle of duality in geometry using the implicit representation enables to "discover" some new formulations or even new theorems. 
    Basic geometric entities and operators are presented by TAB.\ref{tab:my_label_1} and TAB.\ref{tab:my_label_2}.
    \begin{table}[h!]
        \centering
        \caption{Duality of geometric entities}
        \begin{tabular}{|c|c|c|c|c|c|c|}
        \hline
            \multicolumn{7}{|c|}{Duality of geometric entities}  \\ \hline
            Point in $E^2$ & $\xLeftrightarrow[\text{DUAL}]{\text{}}$ & Line in $E^3$ & \quad \quad &  Point in $E^3$ & $\xLeftrightarrow[\text{DUAL}]{\text{}}$ & Plane in $E^3$ \\ \hline
        \end{tabular}
        \label{tab:my_label_1}
    \end{table}
    \begin{table}[h!]
    \centering
        \caption{Duality of operators}
        \begin{tabular}{|c|c|c|}
            \hline
            \multicolumn{3}{|c|}{Duality of operators}  \\ 
            \hline
            Union  $\cup$ & $\xLeftrightarrow[\text{DUAL}]{\text{}}$ & Intersection  $\cap$ \\ \hline
        \end{tabular}
        \label{tab:my_label_2}
    \end{table}
    It means that intersection computation of two line $p_1$ and $p_2$ in $E^2$ is dual to the computation of a line $p$ given by two points $x_1$ and $x_2$ in  $E^2$ using the homogeneous coordinates. 
    In the $E^2$ case, a point $(X,Y)$ is given by homogeneous coordinates $\mathbf{x}=[x,y:w]^T$ 
    and a line $p: aX+bY+c=0, i.e. ax+by+cw = 0 $ by coefficients $[a,b:c]^T$.
    \\
    The usual solutions lead to:
    \begin{itemize}
        \item $\mathbf{Ax=b}$ in the first case, while
        \item $\mathbf{Ax=0}$ in the second case, as the parameters $a,b:c$ of a line are to be determined.
    \end{itemize}
    It is strange, as these problems are dual problems, but formal descriptions are different, but solved differently.
        \\
    However, if the projective formulation  is used, the both cases are solved as the homogeneous system of linear equations, i.e,  $\mathbf{Ax=0}$.
    Similarly, in the $P^3$ case, i.e.  the computation of the intersection point of three planes is dual to the  computation of a plane given by three points, it leads to a system $\mathbf{Ax=0}$, too.

\section{Solution of linear systems of equations}
% -----------------------------

    The linear system of equations $\mathbf{Ax=b}$ can be transformed to the homogeneous system of linear equations, i.e. to the form $\mathbf{D\xi=0}$, where  $\mathbf{D}=[\mathbf{A} |\mathbf{-b}]$,  $\mathbf{\xi}=[\xi_1,...,\xi_n : \xi_w ]^T$, $x_i$ = $\xi_i$ / $\xi_w$, $i=1,...,n$. If $\xi_w \mapsto 0$ then the solution is in infinity and the vector $(\xi_1,...,\xi_n)$ gives the "direction", only.
    
    As the solution of a linear system of equations is equivalent to the outer product (generalized cross-vector) of vectors formed by rows of the matrix $\mathbf{D}$, the solution of the system $\mathbf{D\xi=0}$ is defined as:
    \begin{equation} \label{eq:ksi-1}
        \mathbf{\xi} = \mathbf{d}_1 \wedge \mathbf{d}_2 \wedge ... \wedge \mathbf{d}_n \hspace{2cm} 
        \mathbf{D } \mathbf{\xi} = 0 
        \hspace{1cm} ,\ i.e.  \hspace{1cm}
        \mathbf{[ A |-b] } \mathbf{\xi} = 0 
    \end{equation}
    where: $\mathbf{d}_i$ is the $i$-th row of the matrix $\mathbf{D}$, i.e.  $\mathbf{d}_i=(a_{i1},...,a_{in},-b_i)$, $i=1,...,n$. 
    The application of the projective extension of the Euclidean space enables us to transform the non-homogeneous system of linear equations $\mathbf{Ax=b}$ to the homogeneous linear system $\mathbf{D\xi=0}$, i.e.:
    \begin{equation}
      \begin{bmatrix}
        a_{11} & \dotsb & a_{1n}\\
        \vdots & \ddots & \vdots\\
        a_{n1} & \dotsb & a_{nn}\\ 
        \end{bmatrix} 
        \begin{bmatrix}
            x_{1} \\ \vdots \\     x_{n} \\ 
        \end{bmatrix}
         =
        \begin{bmatrix}
            b_{1} \\ \vdots \\     b_{n} \\ 
        \end{bmatrix}
        % -----
        \quad \xLeftrightarrow[\text{conversion}] \quad \quad
        % ----
        \begin{bmatrix}
        a_{11} & \dotsb & a_{1n} & -b_1\\
        \vdots & \ddots & \vdots & \vdots\\
        a_{n1} & \dotsb & a_{nn} & -b_n\\ 
        \end{bmatrix} 
        \begin{bmatrix}
            \xi_1 \\ \vdots \\     \xi_{n} \\ \xi_{w} \\
        \end{bmatrix}
         =
        \begin{bmatrix}
            0 \\ \vdots \\     0 \\ 
        \end{bmatrix}
    \end{equation}
    It should be noted, that the row rank of the matrix $\mathbf{A}$, $n \times (n+1)$, in the $\mathbf{Ax=0}$ case must be $n$.
    \\ \\
    There are the following important results:
    \begin{itemize}
        \item a solution of a linear system of equations is formally the same for both types, i.e. homogeneous linear systems $\mathbf{Ax=0}$ and non-homogeneous systems $\mathbf{Ax=b}$,
        \item a solution of a linear system of equations is given in the \textit{analytical} form as 
        \[ \mathbf{\xi} = \mathbf{d}_1 \wedge \mathbf{d}_2 \wedge ... \wedge \mathbf{d}_n  \] 
        and relevant operations known for vector space can be used for future processing without numerical evaluation of the linear system.
    \end{itemize}
%

% ===========================
\section{New Geometric Transformation}
    General linear transformations are more complex, especially if the \textit{dot product} and \textit{outer-product} (equivalent to the \textit{cross-product} or \textit{skew-product} 
    in the $E^3$ case) are used.
% -------------------------
    \\
    \\
    \textbf{Basic rules}   \\   
    In the case of the cross-product in the $E^3$ space, the following identity is valid, see Wiki\cite{enwiki:Cross-product}:
    \begin{equation} \label{Eq:cross-product-E3}
        (\mathbf{M} \mathbf{a} ) \times (\mathbf{M} \mathbf{b} ) =
        det(\mathbf{M}) (\mathbf{M}^{-1})^T (\mathbf{a} \times \mathbf{b})
    \end{equation}
    If the matrix $\mathbf{M}$ is orthonormal, then $det(\mathbf{M})=1$,
    the transformation in the Eq.\ref{Eq:cross-product-E3} can be  simplified to:
    \begin{equation}
        (\mathbf{Q} \mathbf{a} ) \times (\mathbf{Q} \mathbf{b} ) =
        \mathbf{Q} ~ (\mathbf{a} \times \mathbf{b})
    \end{equation}
    However, for the $n$-dimensional space and the \textit{outer-product} applications, more general rules can be derived:
    \begin{equation} \label{Eq:cross-product-En}
    \begin{split}
        (\mathbf{M} \mathbf{a} ) \wedge (\mathbf{M} \mathbf{a}_2 ) \wedge  
        \ldots \wedge (\mathbf{M} \mathbf{a}_n)
        = \\
        det(\mathbf{M})^{n-1} (\mathbf{M}^{-1})^T (\mathbf{a}_1 \wedge \mathbf{a}_2 \wedge \ldots \wedge \mathbf{a}_n)
    \end{split}
    \end{equation}
    The presented rules are important as  they enable to handle geometric transformations with lines, planes and  normal vectors. 
    It should be noted that the normal vector of a plane or triangle is actually a bivector and geometric transformation have to respect Eq.\ref{Eq:cross-product-En}.  
    \\
    \\
        Note, that the row vectors $\mathbf{r}_i$, resp. $\mathbf{s}_i$, are the $i$-th row of the matrix 
    $\mathbf{R}$, resp. $\mathbf{S}$. \\
    Then the result of the \textit{geometric product} can be represented as:
    \begin{equation}
        \begin{split}
            \mathbf{a}\mathbf{b} = 
            \mathbf{a} \cdot \mathbf{b} + \mathbf{a} \wedge \mathbf{b} 
            \quad , \quad
            \mathbf{a}\mathbf{b} \Leftrightarrow 
            \mathbf{a} \otimes \mathbf{b} =
            \mathbf{a} ~\mathbf{Q} ~\mathbf{b} 
        \end{split}
    \end{equation}
    where the matrix $\mathbf{Q} = \{ q_{ij} \}$, $i,j=1,\ldots,n$, 
    $q_{ij} = \mathbf{e}_i\mathbf{e}_j$ and
    $\{\mathbf{e}_i\}_{i=1}^n$ are the orthonormal basis vectors in $E^n$.
    \\
    If different transformations $\mathbf{R}$ and $\mathbf{S}$ are applied on the vectors
    $\mathbf{a}$ and $\mathbf{b}$, then:
    \begin{equation}
        \begin{split}
            (\mathbf{R}\mathbf{a}) (\mathbf{S}\mathbf{b}) = 
            (\mathbf{R}\mathbf{a}) \cdot (\mathbf{S}\mathbf{b}) + 
                (\mathbf{R}\mathbf{a}) \wedge (\mathbf{S}\mathbf{b})
            \quad \quad \quad
            (\mathbf{R}\mathbf{a}) (\mathbf{S}\mathbf{b}) \Leftrightarrow 
            (\mathbf{R}\mathbf{a}) \otimes (\mathbf{S}\mathbf{b}) =
            \mathbf{a} ~\mathbf{W} ~\mathbf{b} 
        \end{split}
    \end{equation}
    where the matrix $\mathbf{W} = \{ w_{ij} \}$, $i,j=1,\ldots,n$, $\mathbf{Q}$ is a matrix containing the basis vectors, too.
    \\ 
    Note, that the elements $w_{ij}$ of the matrix $\mathbf{W}$ are given as $\mathbf{e}_i  \mathbf{e}_j$ not shown explicitly.
    \begin{equation}
    \begin{split}
        (\mathbf{R}\mathbf{a}) \otimes (\mathbf{S}\mathbf{b}) = \mathbf{Q}  = 
            \mathbf{a} ~\mathbf{W} ~\mathbf{b} 
        \quad \quad \quad
        q_{ij} = (\mathbf{r_i} \cdot \mathbf{a}) ~ (\mathbf{s_j} \cdot \mathbf{b}) = 
                \mathbf{a} ~ ( \mathbf{r_i} \otimes \mathbf{s_j} ) ~ \mathbf{b}
        \quad \quad \quad
        \mbox{and}\quad w_{ij} =  ( \mathbf{r_i} \otimes \mathbf{s_j} )
    \end{split}
    \end{equation}
    where $ w_{ij} =\mathbf{r_i} \otimes \mathbf{s_j} $, $i,j=1,\dots,n$.
    \\
    Using the dual algebraic adjustments using the multilinearity property WiKi \cite{enwiki:Multilinear-2021} the dual formulation is formed as:
    \begin{equation} \label{Eq:Transform-2}
        \begin{split}
            (\mathbf{R}\mathbf{a}) \otimes (\mathbf{S}\mathbf{b}) = 
                (\mathbf{r_i} \cdot \mathbf{a}) ~ (\mathbf{s_j} \cdot \mathbf{b}) =
                \mathbf{r_i} ~ (\mathbf{a} \otimes\mathbf{b}) ~ \mathbf{s_j} 
                \\
            (\mathbf{R}\mathbf{a}) \otimes (\mathbf{S}\mathbf{b}) = \mathbf{Q}  = 
            \mathbf{r} ~\mathbf{W} ~\mathbf{s}
            \quad \quad \quad
        q_{ij} = (\mathbf{r_i} \cdot \mathbf{a}) ~ (\mathbf{s_j} \cdot \mathbf{b}) = 
                \mathbf{r_i}  ~ (  \mathbf{a} \otimes \mathbf{b} ) ~   \mathbf{s_j} 
            \quad \quad \quad
        \mbox{and}\quad w_{ij} =  ( \mathbf{a} \otimes \mathbf{b} )
        \end{split}
    \end{equation}
    It should be noted that the diagonal of the matrix 
    $\mathbf{W}$ contains elements of the \textit{inner product}.
    The non-diagonal elements represent parts of bivectors of the given $n$-dimensional space.
    \\
    This formulation Eq.\ref{Eq:Transform-2} has the advantage that for the given constant transformations $\mathbf{R}$ and $\mathbf{S}$ the matrix is $\mathbf{W}$ is \textit{constant}.

% --------------------
\section{Geometric examples \protect}
    Let us demonstrate the power of the geometric algebra on simple geometric problems and their simple solutions, e.g. intersection of two lines in the $E^2$ space, intersection of three planes in the $E^3$ space and dual problems, intersection of two planes in the $E^3$, barycentric coordinates, etc.
%
% -----------------------------
    \\
    \\
    \textbf{Two and dimensional examples}   \\
    The direct consequence of the principle of duality is that, the intersection point $\mathbf{x}$ of two lines in $E^2$  $\mathbf{p}_1,\mathbf{p}_2$, resp. a line $\mathbf{p}$ passing two given points $\mathbf{x}_1,\mathbf{x}_2$, is given as:
    \begin{equation}
        \mathbf{x}=\mathbf{p}_{1} \wedge \mathbf{p}_2 
        \xLeftrightarrow[\text{DUAL}]{\text{}} 
        \mathbf{p}=\mathbf{x}_1 \wedge \mathbf{x}_2
    \end{equation}
    where $\mathbf{p}_i=[a_i,b_i:c_i ]^T$,  $\mathbf{x}=[x,y:w]^T$ ($w$ is the homogeneous coordinate), $i=1,2$; similarly in the dual case.
    \\
    In the case of the $E^3$ space, a point is dual to a plane and vice versa. It means that the intersection point $\mathbf{x}$ of three planes $\mathbf{\rho}_1$,$\mathbf{\rho}_2$,$\mathbf{\rho}_3$, resp. a plane $\mathbf{\rho}$ passing three given points $\mathbf{x}_1,\mathbf{x}_2,\mathbf{x}_3$ is given as:
    \begin{equation}
        \mathbf{x}=\mathbf{\rho}_{1} \wedge \mathbf{\rho}_2 \wedge \mathbf{\rho}_3 \xLeftrightarrow[\text{DUAL}]{\text{}} 
        \mathbf{\rho}=\mathbf{x}_1 \wedge \mathbf{x}_2 \wedge \mathbf{x}_3
    \end{equation}
    where $\mathbf{x}=[x,y,z:w]^T$, $\mathbf{\rho}_i=[a_i,b_i,c_i:d_i ]^T$, $ i=1,2,3$. 
    \\
    It can be seen that the above formulae are equivalent to the application of the outer product (“extended” cross-product), which in natively supported by GPU architecture. 
    \\
    For an intersection computation, we get:
    \begin{equation}
        \mathbf{x}=\mathbf{p}_1 \wedge \mathbf{p}_2 = 
        \begin{bmatrix}
            \mathbf{e}_1 & \mathbf{e}_2 & \mathbf{e}_w\\
            a_1 & b_1 & c_1 \\
            a_2 & b_2 & c_2 
        \end{bmatrix} 
        \hspace{1cm}
        \mathbf{x}=\mathbf{\rho}_1 \wedge \mathbf{\rho}_2 \wedge \mathbf{\rho}_3 = 
        \begin{bmatrix}
            \mathbf{e}_1 & \mathbf{e}_2 & \mathbf{e}_3 & \mathbf{e}_w\\
            a_1 & b_1 & c_1 & d_1\\
            a_2 & b_2 & c_2 & d_2\\
            a_3 & b_3 & c_3 & d_3
        \end{bmatrix}
        \end{equation}
        Due to the principle of duality, a dual problem solution is given as:
        \begin{equation}
            \mathbf{p} = \mathbf{x}_1 \wedge \mathbf{x}_2 = 
            \begin{bmatrix}
                \mathbf{e}_1 & \mathbf{e}_2 & \mathbf{e}_w\\
                x_1 & y_1 & w_1 \\
                x_2 & y_2 & w_2 
            \end{bmatrix} 
        \hspace{2cm} 
        \mathbf{\rho}=\mathbf{x}_1 \wedge \mathbf{x}_2 \wedge \mathbf{x}_3 = 
        \begin{bmatrix}
            \mathbf{e}_1 & \mathbf{e}_2 & \mathbf{e}_3 & \mathbf{e}_w\\
            x_1 & y_1 & z_1 & w_1\\
            x_1 & y_2 & z_2 & w_2\\
            x_3 & y_3 & z_3 & w_3
        \end{bmatrix}
    \end{equation}
    The above-presented formulae prove the strength of the formal notation of the geometric algebra approach. 
    Therefore, there is a natural question, what is the more convenient computation of the geometric product, as computation with the outer product, i.e. extended cross product, using the basis vector notation approach is not easy.
%    
% -------------------------
    \\
    \\
    \textbf{Barycentric coordinates}   \\
    The barycentric coordinates are often used in many
    applications, not only in geometry. The barycentric coordinates
    computation leads to a solution of a system of linear
    equations. 
    However, a solution of a linear system equations is equivalent to the outer product Skala\cite{Skala2005905}\cite{Skala2006625}. 
    Therefore, it is possible to compute the barycentric
    coordinates using the outer product, which is recommendable especially for the GPU oriented applications.

    Let us consider the $E^2$ case and the barycentric interpolation between three points (vertices) $\mathbf{x}_i=[x_i,y_i:w_i]^T$, $i=1,\ldots,3$, of the given triangle, and vectors:
    \begin{equation}
    \begin{split}
        \mathbf{\xi} = [x_1,x_2,x_3:x]^T \quad
        \mathbf{\eta} = [y_1,y_2,y_3:y]^T \quad
        \mathbf{\omega} = [w_1,w_2,w_3:w]^T 
    \end{split}
    \end{equation}
    Then the barycentric coordinates $\mathbf{\mu}$ in the homogeneous coordinates of the point $\mathbf{x} = [x,y:w]^T$ are given as:
    \begin{equation}
    \begin{split}
        \mathbf{\mu} = \mathbf{\xi} \wedge \mathbf{\eta} \wedge \mathbf{\omega}
    \end{split}
    \end{equation}
    where $\mathbf{\mu} = [\mu_1, \mu_2,\mu_3: \mu_w]^T $
    and the barycentric coordinates in the Euclidean space $\mathbf{\lambda}$ are given as:
    \begin{equation}
    \begin{split}
        \mathbf{\lambda} = (\lambda_1,\lambda_2,\lambda_3) =
        (-\frac{\mu_1}{\mu_w},-\frac{\mu_2}{\mu_w},-\frac{\mu_3}{\mu_w})
    \end{split}
    \end{equation}   
    Similarly, for other dimensions, see Skala\cite{Skala2008120} for details. 
    How simple and elegant solution!
    
    It can be seen that the presented computation of barycentric
    coordinates is simple, convenient for GPU or SSE
    application. 
    Even more, as we have assumed from the very
    beginning, there is no need to convert coordinates of points
    from the homogeneous coordinates to the Euclidean
    coordinates. 
     As a direct consequence of that is that we save
    lot of division operations and also increase robustness of the
    computation.
%
% -------------------------------------
    \\
    \\
    \textbf{Pl\"ucker coordinates}    \\
    There are two other geometric problems in the $E^3$ case, i.e. Intersection of two planes in the $E^3$ space and its dual problem, 
    i.e. a line given by two points in the $E^3$ space.
    
    The geometric product of vectors representing two planes $ \mathbf{\rho}_i=[a_i,b_i,c_i:d_i]^T$, resp. two points  
    $\mathbf{x}_i=[x_i,y_i,z_i:w_i]^T$, $i=1,\ldots,3$, using the homogeneous coordinates is given using the anti-commutative tensor product as:
    \begin{center}
    \begin{tabular}{ |c|c|c|c|c|c|} 
     \hline
        ~$\mathbf{\rho}_1\mathbf{\rho}_2$~ & $a_2$ & $b_2$ & $c_2$ & $d_2$ \\ \hline
        $a_1$ & $a_1 a_2$ & $a_1 b_2$ & $a_1 c_2$ & $a_1 d_2$  \\ \hline
        $b_1$ & $b_1 a_2$ & $b_1 b_2$ & $b_1 c_2$ & $b_1 d_2$  \\ \hline
        $c_1$ & $c_1 a_2$ & $c_1 b_2$ & $c_1 c_2$ & $a_1 d_2$  \\ \hline
        $d_1$ & $d_1 a_2$ & $d_1 b_2$ & $d_1 c_2$ & $d_1 d_2$  \\ \hline
    \end{tabular}
    \hspace{1cm}
    \begin{tabular}{ |c|c|c|c|c|c|} 
     \hline
        ~$\mathbf{x}_1\mathbf{x}_2$~ & $x_2$ & $y_2$ & $z_2$ & $w_2$ \\ \hline
        $x_1$ & $x_1 x_2$ & $x_1 y_2$ & $x_1 z_2$ & $x_1 w_2$  \\ \hline
        $y_1$ & $y_1 x_2$ & $y_1 y_2$ & $y_1 z_2$ & $y_1 w_2$  \\ \hline
        $z_1$ & $z_1 x_2$ & $z_1 y_2$ & $z_1 z_2$ & $x_1 w_2$  \\ \hline
        $w_1$ & $w_1 x_2$ & $w_1 y_2$ & $w_1 z_2$ & $w_1 w_2$  \\ \hline
    \end{tabular}
    \end{center}
    However, the question is how to compute a line $\mathbf{p}\in E^3$ given as an intersection of two planes $\mathbf{\rho}_1$, $\mathbf{\rho}_2$, which is dual to a line determination given by two points $\mathbf{x}_1$, $\mathbf{x}_2$ as those problems are dual. 
    
    The parametric solution can be easily obtained using standard Plücker coordinates, however computation and formula are complex and not easy to understand. 
    \begin{equation}
    \begin{split}
        \mathbf{q}(t)=\frac{\mathbf{\omega} \wedge \mathbf{v}  } {||\mathbf{\omega}||^2}+\mathbf{\omega}\ t \hspace{2cm}
        \mathbf{L}=\mathbf{x}_1 \mathbf{x}_2^T-\mathbf{x}_2 \mathbf{x}_1^T 
        \quad , resp. \quad
        \mathbf{L}=\mathbf{\rho}_1 \mathbf{\rho}_2^T-\mathbf{\rho}_2 \mathbf{\rho}_1^T 
        \\
        \mathbf{\omega} = [l_{41},l_{42},l_{43}]^T 
        \hspace{2cm} 
        \mathbf{v} = [l_{23},l_{31},l_{12}]^T 
    \end{split}
    \end{equation}
    where $l_{ij}$ are the Pl\"ucker coordinates and $\mathbf{q}(t)$ is a line in $E^3$ in the parametric form.  \\
    For the case of intersection of two planes the principle of duality can be applied directly.
    \\
    However, using the geometric algebra, principle of duality and projective representation, we can directly write:
    \begin{equation}
        \mathbf{p} = \mathbf{\rho_1} \wedge \mathbf{\rho_2} \xLeftrightarrow[\text{DUAL}]{\text{}}  
        \mathbf{p} = \mathbf{x}_1 \wedge \mathbf{x}_2
    \end{equation}
    It can be seen that the formula given above keeps the duality in the final formulae, too.
    
    From the formal point of view, the geometric product for both cases is given as:
    \begin{equation}
    \begin{split}
        \mathbf{\rho}_1\mathbf{\rho}_2 \xLeftrightarrow[\text{repr}]{\text{}} 
        \mathbf{\rho}_1 \otimes \mathbf{\rho}_2 = 
    \begin{bmatrix}
        a_1 a_2 & a_1 b_2 & a_1 c_2 & a_1 d_2\\
        b_1 a_2 & b_1 b_2 & b_1 c_2 & b_1 d_2\\
        c_1 a_2 & c_1 b_2 & c_1 c_2 & c_1 d_2\\
        d_1 a_2 & d_1 b_2 & d_1 c_2 & d_1 d_2\\
    \end{bmatrix}
    \xLeftrightarrow[\text{DUAL}]{\text{}}
        \mathbf{x}_1 \mathbf{x}_2 \xLeftrightarrow[\text{repr}]{\text{}} \mathbf{x}_1 \otimes \mathbf{x}_2 = 
    \begin{bmatrix}
        x_1 x_2 & x_1 y_2 & x_1 z_2 & x_1 w_2\\
        y_1 x_2 & y_1 y_2 & y_1 z_2 & y_1 w_2\\
        z_1 x_2 & z_1 y_2 & z_1 z_2 & z_1 w_2\\
        w_1 x_2 & w_1 y_2 & w_1 z_2 & w_1 w_2\\
    \end{bmatrix}
    \end{split}
    \end{equation}
    It means that we have the computation of the Plücker coordinates for the both cases, i.e. for computation of a line $\mathbf{p}=\mathbf{\rho}_1 \wedge \mathbf{\rho}_2$ or $\mathbf{p}=\mathbf{x}_1 \wedge \mathbf{x}_2$ is given as a union of two points in $E^3$ and as an intersection of two planes in $E^3$ using the projective representation and the principle of duality. 
    It should be noted that the given approach offers significant simplification of computation of the Plücker coordinates as it is simple and easy to derive and explain.
    It also uses vector-vector operations, which is especially convenient for SSE and GPU application one code sequence for the both cases.
    
    As the Plücker coordinates are also in mechanical engineering applications, especially in robotics due to its simple displacement and momentum specifications, and in other fields simple explanation and derivation is another very important argument for GA approach application.
\newpage

% -----------------------------
\subsection{Conditionality - Angular criterion}
% -----------------------------
    The robustness and reliability of a solution of linear systems of equations $\mathbf{Ax=b}$ is a critical issue and approaches are mostly based on eigenvalues of a matrix evaluation.
    In our case, the both types of the linear systems of equations, i.e. $\mathbf{Ax=b}$ ($\mathbf{A}$ is $n \times n$) and $\mathbf{Ax=0}$ ($\mathbf{A}$ is $n \times (n+1)$), actually have the same form $\mathbf{Ax=0}$ ($\mathbf{A}$ is $n \times (n+1)$, if the projective representation is used.
    
    Let us consider the $i$ row $\mathbf{a}_i$ of the matrix $\mathbf{A}$ and the $i$ row $\mathbf{\alpha}_i$ of the matrix $[\mathbf{A}|\mathbf{-b}]$ in the case of the linear system $\mathbf{Ax=b}$.

\begin{figure}[!ht]
        \centering
        \includegraphics[width=6cm]{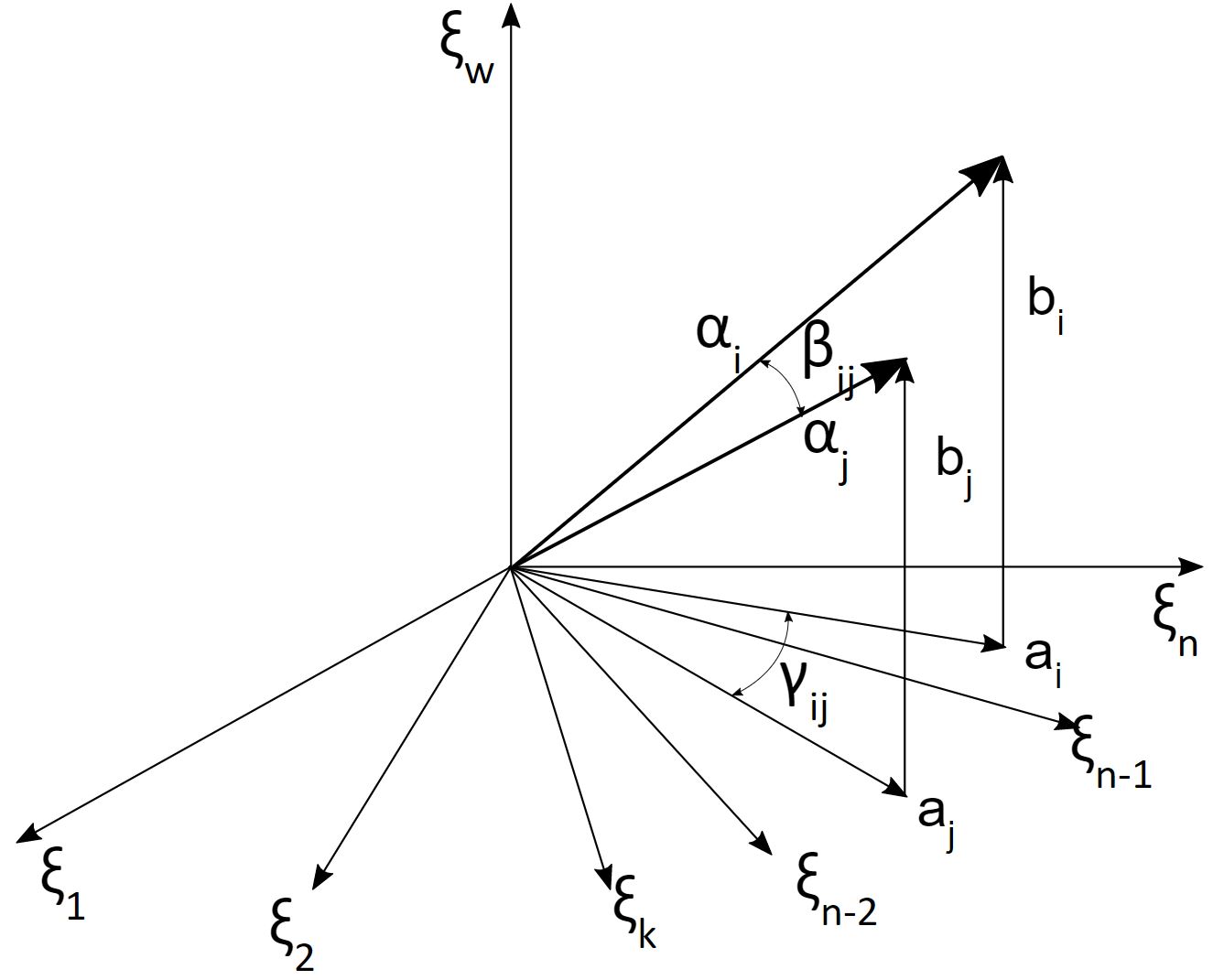}
        \caption{Difference between matrix and linear system conditionality}
        \label{fig:Fig-3} \label{coditionality-geom}
    \end{figure}

    The solution $\mathbf{\tau}$ of the linear system is given as   
    $\mathbf{\tau} = \mathbf{\alpha}_1 \wedge \mathbf{\alpha}_2 \wedge ... \wedge \mathbf{\alpha}_n$ and
    the vector $\mathbf{\xi} = [\xi_1,\ldots,\xi_n:\xi_w]^T $ represents the orthonormal basis, see Fig.\ref{coditionality-geom}.
    The angle $\gamma_{ij}$ is the angle between two vectors $\mathbf{a}_i$ and $\mathbf{a}_j$. 
    If there, in the case of 4$\mathbf{Ax=b}$,  exists $i,j$ so that the value of $\gamma_{ij}$ is close, resp. equal to zero, the matrix $\mathbf{A}$ is close to singular, resp. is singular.
    \\
    In the case of the homogeneous linear system, i.e. $\mathbf{Ax=0}$ the rows $\mathbf{\alpha}_i$ and $\mathbf{\alpha}_j$ of the matrix $\mathbf{A}$ are close to, resp. linearly dependent, the 
    $\beta{ij}$ is close, resp. close to zero, see Fig.\ref{coditionality-geom}.
    
    Therefore, it is possible to see the differences between the matrix conditionality and conditionality (solvability) of a linear system of equations, see Fig.\ref{fig:Fig-3}. 
    
    The eigenvalues are usually used  and the ratio $rat_{\lambda}=|\lambda_{max}| /  |\lambda_{min}|$ \& $\lambda_i \in C^1$ is mostly used as a criterion of the  conditionality. 
    If the ratio $rat_{\lambda}$ is high, the matrix is said to be ill-conditioned.
    It should be noted that the computation of eigenvalues is approx. of the $O(n^3)$ time complexity, i.e. extremely slow especially in the case of large data with a large span of data. 
    The second approach is based on an error 
    $r = \|\mathbf(A) \| \|\mathbf(A^{-1}) \|$ computation.
    \\
    There are two cases, which are needed to be taken into consideration:
    % -------------------
    \begin{itemize}  % first level
        \item non-homogeneous systems of linear equations, i.e. $\mathbf{Ax=b}$. In this case, the matrix conditionality is considered as a criterion for the solvability of the linear system of equations. It depends on the matrix $\mathbf{A}$ properties, i.e. on the eigenvalues given as $\det(\mathbf{A} - \lambda \mathbf{I})=0  $. 
        
        A conditionality number $\kappa(\mathbf{A})=|\lambda_{max}| /  |\lambda_{min}| $ is usually used as the solvability criterion.
        Let us consider a simple example:
        \begin{equation} \label{eq:diagonal-1}
            \begin{bmatrix}
            10^2 & 0 & 0  \\
             0  & 10^0 & 0  \\
             0 & 0 & 10^{-2}\\ 
            \end{bmatrix} 
            \begin{bmatrix}
                x_{1} \\ \vdots \\     x_{3} \\ 
            \end{bmatrix}
             =
            \begin{bmatrix}
                b_{1} \\ \vdots \\     b_{3} \\ 
            \end{bmatrix}
            \quad 
            \begin{bmatrix}
             0       & 0    & 10^2  \\
             0       & 10^0 & 0  \\
             10^{-2} & 0    & 0\\ 
            \end{bmatrix} 
            \begin{bmatrix}
                \xi_{1} \\ \vdots \\     \xi_{3} \\ 
            \end{bmatrix}
             =
            \begin{bmatrix}
                \beta_{1} \\ \vdots \\     \beta_{3} \\ 
            \end{bmatrix}
            \quad 
            \xi_i = x_{i-4}, 
            \beta_i = b_{i-4} \quad i=1,2,3
        \end{equation}
        In the case of the Eq.\ref{eq:diagonal-1}, the matrix conditionality is $\kappa(\mathbf{A})=10^2/10^{-2}=10^4$. 
        However, if the $1^{st}$ row is multiplied by $10^{-2}$ and the $3^{rd}$ row is multiplied by $10^2$, then the conditionality is $\kappa(\mathbf{A})=1$. 
        \\
        It can be seen that both linear systems, see Eq.\ref{eq:diagonal-1}, represent the equivalent problem.
        % ----
        \item the homogeneous system of equations $\mathbf{Ax=0}$, when the system of linear equations  $\mathbf{Ax=b}$ is expressed in the projective space. In this case, the vector $\mathbf{b}$ is taken into account and the bivector area and the bivector angles properties can be used for solvability evaluation.
    \end{itemize} % first level
    The angular conditionality can be express as 
    \begin{equation}
    \begin{split}
        \kappa_\gamma(\mathbf{A}) = \frac{ min{~\gamma_{ij}} }{ max{~\gamma_{kl}} } 
        \quad \quad
        \kappa_\gamma(\mathbf{A}) \in <0,1>
        \quad \quad
        \gamma_{ij} = \arccos \frac{|\mathbf{a}_i \cdot \mathbf{a}_j|}{\| \mathbf{a}_i \| \| \mathbf{a}_j \|}  \\  %
        \kappa_\beta([\mathbf{A} | \mathbf{-b}]) = \frac{ min{~\beta_{ij}} }{ max{~\beta_{kl}} } 
        \quad \quad
        \kappa_\beta([\mathbf{A} | \mathbf{-b}]) \in <0,1>
        \quad \quad
        \beta_{ij} = \arccos
        \frac{\mathbf{|\alpha}_i \cdot \mathbf{\alpha}_j|}{\| \mathbf{\alpha}_i \| \| \mathbf{\alpha}_j \|}
        \\ i,j,k,l = 1,\ldots,n ~ \& ~ i \ne j ~\&~ k \ne l  
    \end{split}
    \end{equation}
    \quad
    \\
    The proposed \textit{angular} conditionality criterion is invariant to the row multiplications, while only the column multiplication (representing a change of the physical units of the $x_i$) changes the angles of the bivectors. \\
    \\
    There are several significant consequences:
    % ------------------
    \begin{itemize}
        \item the solvability of a linear system of equations can be improved by the column multiplications, only if unlimited precision is considered.
        Therefore, the matrix-based pre-conditioners might not solve the solvability problems and might introduce additional numerical problems, see Chen\cite{Chen2005}, Benzi\cite{Benzi2002},
        \item the precision of computation is significantly influenced by addition and subtraction operations in the floating-point representation\cite{wiki:IEEE-754}\cite{IEEE754-219}, as the exponents must be the same for those operations with mantissa. Also, the multiplication and division operations using exponent change by $2^{\pm k}$ should be preferred.
    \end{itemize}
    It should be noted that, the $\log_{2}(*)$, resp. $\log_{10}(*)$  function is to be used for the practical use as the exponent value is interesting only for the conditionality assessment.
%
% -----------------------------------
\section{Conclusion}
    This contribution presents a short introduction to geometric algebra principles, \textit{geometric product}, \textit{outer product},  \textit{inner product} and the anti-commutative tensor product for efficient computation. 
    It presents the following main contributions:
    \begin{itemize}
        \item the geometric product and the outer product extension for applications in the projective space,
        \item the equivalence of a solution of linear systems of equations with the  \textit{outer product} application,
        \item a solution of the linear system of equations $\mathbf{Ax=b}$, resp. $\mathbf{Ax=0}$ is available in the \textit{analytical} form, i.e. \\ 
        $\mathbf{\xi} = \mathbf{d}_1 \wedge \mathbf{d}_2 \wedge ... \wedge \mathbf{d}_n$
        and relevant operations known for vector space can be used for future processing without numerical evaluation the linear system,
        \item a unique approach to a solution of homogeneous systems linear of equations, i.e. $\mathbf{Ax=0}$, and non-homogeneous  systems linear of equations, i.e. $\mathbf{Ax=b}$ using the outer product,
        \item an application of the principle of duality in solving selected geometrical problems, e.g. computation of the barycentric coordinates, simplification of the Pl\"ucker coordinates, etc.,
        \item a new approach to the evaluation of a matrix conditionality based on angular ratios of row vectors of the given extended matrix $[~\mathbf{A} | \mathbf{-b}~]$ in the case of $\mathbf{A}\mathbf{x}=\mathbf{b}$, 
        resp. of a matrix $\mathbf{A}$ in the case of 
        $\mathbf{A}\mathbf{x}=\mathbf{0}$.
    \end{itemize}
%
%-------------------------
\begin{acknowledgments}
    The author would like to thank to colleagues and students at the University of West Bohemia (Czech Republic), Shandong University and Zhejiang University (China) for their critical comments and constructive suggestions, and to anonymous reviewers for their valuable comments and hints provided.
\end{acknowledgments}

%\nocite{*}
% \bibliographystyle{apsrev4-1} 
%\bibliographystyle{unsrt}
\bibliography{outpub.bbl}
%  \bibliography{mybibfile}% Produces the bibliography via BibTeX.

\section*{Appendix}
    The GPU implementation of the outer product for the $E^3$ case using the homogeneous coordinate is quite simple. 
    It should be noted that only 4 clocks for the \textit{outer product} and 4 clocks for the \textit{inner product} are needed.
    \begin{verbatim}
        float4 a;
        a.x = dot(x1.yzw, cross(x2.yzw, x3.yzw));
        a.y = - dot(x1.xzw, cross(x2.xzw, x3.xzw));
        a.z = dot(x1.xyw, cross(x2.xyw, x3.xyw));
        a.w = - dot(x1.xyz, cross(x2.xyz, x3.xyz));
        return a;  
    \end{verbatim}
    or more compactly as:
    \begin{verbatim}
        float4 cross_4D(float4 x1, float4 x2, float4 x3)
        return(
            dot(x1.yzw, cross(x2.yzw, x3.yzw)),   
            - dot(x1.xzw, cross(x2.xzw, x3.xzw)),   
            dot(x1.xyw, cross(x2.xyw, x3.xyw)),    
            - dot(x1.xyz, cross(x2.xyz, x3.xyz)) 
        );     
    \end{verbatim}

\end{document}